\documentclass[12pt]{article}
\usepackage{fancyhdr, array,calc,graphicx,url,tabularx}
\usepackage{geometry}
\usepackage{latexsym}
\usepackage{amssymb}
\usepackage{amsmath}
\usepackage{makeidx}
\usepackage{enumerate}
\usepackage{verbatim}
\usepackage{tikz}
\usepackage[colorlinks=true,linkcolor=blue]{hyperref}

\usepackage{changepage}

\usepackage[utf8]{inputenc}
\usepackage{tikz}
\tikzset{
     block/.style={rectangle, draw, fill=red!40, text width=6em,
                   text centered, rounded corners, minimum height=3em},
     arrow/.style={-{Stealth[]}}
     }

\makeindex

{
   \end{minipage}
   \vspace*{\stretch{3}}
   \clearpage
}

\def\undertilde#1{{\baselineskip=0pt\vtop
  {\hbox{$#1$}\hbox{$\scriptscriptstyle\sim$}}}{}}

\newcommand{\utilde}{\undertilde}

\renewcommand{\gg}{\gamma}

\newcommand{\bR}{{\mathbb{R}}}

\newcommand{\rest}{\restriction}

\newcommand{\forces}{\Vdash}

\renewcommand{\models}{\vDash}
\newcommand{\powerset}{{\wp}}
\newcommand{\cp}{{\rm crit }}

\newcommand{\cf}{{\rm cf}}

\newtheorem{theorem}{Theorem}[section]
\newtheorem{proposition}[theorem]{Proposition}
\newtheorem{definition}[theorem]{Definition}

\newtheorem{lemma}[theorem]{Lemma}
\newtheorem{corollary}[theorem]{Corollary}

\newtheorem{sublemma}[theorem]{Sublemma}

\newtheorem{notation}[theorem]{Notation}

\numberwithin{figure}{section}

\newenvironment{proof}{{\it{
Proof.}}}{\nopagebreak\mbox{}{\hfill$\square$}
\par\bigskip}
\newenvironment{Claim}{\noindent{\it
Claim} {}}{}
\newenvironment{Case}{\noindent{\it
Case}}{}

\newcommand{\rslem}[1]{Sublemma~\ref{#1}}

\newcommand{\rthm}[1]{Theorem~\ref{#1}}
\newcommand{\rlem}[1]{Lemma~\ref{#1}}

\newcommand{\rcor}[1]{Corollary~\ref{#1}}
\newcommand{\rdef}[1]{Definition~\ref{#1}}

\def\k{\kappa}
\def\a{\alpha}
\def\b{\beta}
\def\d{\delta}

\def\l{\lambda}

\def\P{{\mathcal{P} }}

\def\H{{\rm{HOD}}}
\def\M{{\mathcal{M}}}
\def\N{{\mathcal{N}}}

\def\T {{\mathcal{T}}}
\def\U{{\mathcal{U}}}

\def\iff{\mathrel{\leftrightarrow}}

\def\and{\mathrel{\kern1pt\&\kern1pt}}

\def\insegeq{\trianglelefteq}

\def\<#1>{\langle\,#1\,\rangle}

 \input xy
 \xyoption{all}
       
\title{$AD_{\mathbb{R}}$ implies that all sets of reals are $\Theta$ universally Baire}

\author{Grigor Sargsyan}

\date{\today}

\begin{document}
\maketitle

The purpose of this note is to establish the following useful theorem. 

\begin{theorem}\label{main theorem} Assume $AD_{\mathbb{R}}$. Then all set of reals are $\Theta$ universally Baire. 
\end{theorem}

We are interested in this theorem for two reasons. It is used in the author's construction  of \textit{generic generators} (see the forthcoming \cite{GG}). Secondly the authors of \cite{AllSetsUB} introduced a derived model construction producing a model of determinacy in which all sets of reals are universally Baire. The above theorem shows how to build such models from determinacy rather than large cardinals\footnote{However, the large cardinal assumption used in \cite{AllSetsUB} is probably optimal.} (see \rcor{a model of all sets are ub}). Along the way of proving \rthm{main theorem}, we fill in some gaps in descriptive inner model theory. In particular, we outline the process of getting $\N^*_x$-like models used in \cite{ATHM} (see \cite[Theorem 2.25]{ATHM}). We still leave one gap. The proof of \rthm{coarse capturing} is not given as the full proof will go beyond the scope of this paper. However, we believe we have given enough outline so that the interested reader can prove it. 

I would like to thank the anonymous referee for a long list of corrections. The author was supported by the NSF Career Award DMS-1352034.

\section{$\Gamma$-Woodins}
With the exception of \rthm{coarse capturing}, all of the material that we review in this section has been developed by the Cabal group, and among other places, can be found in \cite{MouseSetWoodin} (especially \cite[Chapter 3]{MouseSetWoodin}), the Appendix of \cite{ATHM}, \cite{SchSteel} and \cite{CMI}. 

Following \cite[Chapter 3]{MouseSetWoodin}), we say that $\Gamma$ is a \textit{good pointclass}\footnote{$\bR$ is the Baire space.} if
\begin{enumerate}
\item $\Gamma$ is closed under recursive substitution and number quantification,
\item $\Gamma$ is $\omega$-parametrized\footnote{This means that there is a set $U^k\subseteq \omega\times \bR^k$ such that $U^k\in \Gamma$ and for every set $A\in \Gamma$, if $A\subseteq \bR^k$ then there is an integer $n$ such that $x\in A\iff (n, x)\in U^k$.},
\item $\Gamma$ has the scale property,
\item $\Gamma$ is closed under $\exists^\bR$.
\end{enumerate}
$\Sigma^1_2$ is the prototypical good pointclass. Every good pointclass has its boldface version  $\utilde{\Gamma}$ which consists of those sets $A\subseteq \bR^k$ such that for some $U\in \Gamma$ and $x\in \bR$, $U\subseteq \bR^{k+1}$ and $A=U_x=\{\vec{y}:(x, y)\in U\}$.  One can also define the pointclass $\Gamma(x)$ for $x\in \bR$ in a similar manner. 

Each good pointclass has its associated $C_\Gamma$ operator. Here for $x\in \bR$, 
\begin{center}
$C_\Gamma(x)=\{ y\in \bR: y$ is $\Gamma$-definable from $x$ and a countable ordinal$\}$.
\end{center}
The $C_\Gamma$ operator can be extended to sets in $HC$ via the category quantifier. This is done, for instance, in \cite[Chapter 3]{MouseSetWoodin}. \\\\
Below we explain how this is done. Given $a\in HC$ and $f:\omega\rightarrow a$ a surjection, we let $a_f=\{2^n3^m: (n, m)\in \mathbb{N}^2 \wedge f(n)\in f(m)\}$. Given $b\subseteq a$, we let $b_f=\{ n\in \mathbb{N}: f(n)\in b\}$. We can then easily decode $b$ from $(a_f, b_f)$. Indeed, letting $c=\{ (n, m)\in \mathbb{N}^2: 2^n3^m\in a_f\}$, $(a, \in)$ is the transitive collapse of $(\omega, c)$ and $(a, b, \in)$ is the transitive collapse of $(\omega, b_f, c)$. We then set $C_\Gamma(a)$ to be the set of those $b\subseteq a$ such that for co-meager many $f:\omega\rightarrow a$ there is $b_f\in C_\Gamma(a_f)$.
The Harrington-Kechris Theorem gives a nice description of the set $C_\Gamma(a)$.

\begin{theorem}[Harrington-Kechris, \cite{HarKech}]\label{capturing by trees} Assume $AD$. Suppose $T$ is a tree  of a $\Gamma$-scale on a universal $\Gamma$-set. Then for any transitive $a\in HC$,
\begin{enumerate}
\item $C_\Gamma(a)=\{b\subseteq a:$ for all $f:\omega\rightarrow a$, $b_f\in C_\Gamma(a_f)\}$, and
\item $C_\Gamma(a)=\powerset(a)\cap L(T\cup\{a\}, a)$.
\end{enumerate}
\end{theorem}

Next we introduce $\Gamma$-Woodins.
\begin{notation} Below we will use $C_\Gamma^\a$ for the $\a$th iterate of $C_\Gamma$. Thus, $C^2_\Gamma(a)=C_\Gamma(C_\Gamma(a))$. We only need this notion for $\a\leq \omega$. We then set $C^\omega_\Gamma(a)=\cup_{n<\omega}C^n_{\Gamma}(a)$. 

Suppose $T$ is the tree of a $\Gamma$-scale. For each $\a<\omega_1$, we let $\k_\a$ be the $\a$th- infinite cardinal of $L[T,a]$. We can then simply set $C^\a_\Gamma=H_{\k_\a}^{L[T, a]}$. Then, using this definition, we have $C_\Gamma(a)=C^1_\Gamma(a)$.
\end{notation}

\begin{definition}\label{gamma woodin} Given a transitive $P$ model of $ZFC-Replacement$, we say $P$ is a $\Gamma$-Woodin if for some $\d$,
\begin{enumerate}
\item $P\models ``\d$ is the only Woodin cardinal".
\item $P=C^\omega_\Gamma(P)$,
\item for every $P$-inaccessible cardinal $\eta<\d$, $C_\Gamma(V^P_\eta)\models ``\eta$ is not a Woodin cardinal".
\end{enumerate}
We let $\d^P$ be the Woodin cardinal of $P$.
\end{definition}

To make matters simple, we will work with prototypical good pointclass\footnote{The point of this move is to just avoid discussing $Env(\Gamma)$. Readers familiar with this notion do not have to make this move.}. These are pointclass that correspond to $\Sigma^2_1$ of an initial segment of the universe. We say $(A_n: n<\omega)\subseteq \bR^{\omega}$ is a self-justyfying-system (sjs) if  for each $n\in \omega$, 
\begin{enumerate}
\item there is a sequence $(A_{m_k}: k\in \omega)$ that codes a scale on $A_n$,
\item there is $m<\omega$ such that $A_n^c=A_m$.
\end{enumerate} 

Let $T_0$ be the theory 
\begin{enumerate}
\item $ZF-Powerset$,
\item $\Theta$ exists\footnote{More precisely, ``there is an ordinal which is not the surjective image of $\bR$".}, and
\item $V=L_{\Theta^{+}}(C, \bR)$ for some $C\subseteq \bR$.
\end{enumerate}

\begin{definition}\label{very good pointclass} Suppose $\Gamma$ is a good pointclass. We say $\Gamma$ is a very good pointclass if there is a sjs $\vec{A}=(A_n: n\in \omega)$, $\gg<\Theta^{L(\vec{A}, \bR)}$, a $\Sigma_1$-formula $\phi$ and a real $x$ such that $L_\gg(\vec{A}, \bR)$ is the least initial segment of $L(\vec{A}, \bR)$ that satisfies $T_0+\phi(x)$ and $\Gamma=(\Sigma^2_1(\vec{A}))^{L_{\gg}(\vec{A}, \bR))}$. 
\end{definition}

Recall that Martin showed (for instance, see \cite{HomMosc} or \cite{ADPlusBook}) that under $AD_{\mathbb{R}}$, every set of reals is Suslin. It follows that under $AD_{\mathbb{R}}$, for every set $B$ there is a sjs $\vec{A}$ such that $B=A_0$. 

If $\gg$ is as in \rdef{very good pointclass} then $L_{\gg+2}(\vec{A}, \bR)$ has a set of reals not in $L_{\gg+1}(\vec{A}, \bR)$, namely the set of reals coding $L_\gamma(\vec{A}, \bR)$. Thus, $\gg$ ends a weak gap relative to $\vec{A}$. The reason for working with a sjs rather than any set $A$ is that if $\vec{A}$ is a sjs then $L(\vec{A}, \bR)$ has the same properties as $L(\bR)$. In particular, it admits a scales analysis (see \cite{ScalesLR}).

From now on we will write \textit{vg} for very good. Suppose $\Gamma$ is a $vg$-pointclass. We then let $M_\Gamma=L_\gg(\vec{A}, \bR)$ be as in \rdef{very good pointclass} witnessing very goodness of $\Gamma$. We say that $M_\Gamma$ is the \textit{parent} of $\Gamma$. Clearly, $M_\Gamma$ has many representations as $L_\gg(\vec{A}, \bR)$. However, the resulting ambiguity is not problematic and we will not dwell on it. 

Notice that we have that
\begin{proposition} if $\Gamma$ is vg-pointclass and $M_\Gamma=L_\gg(\vec{A}, \bR)$ is its parent then for any countable transitive $a$, $C_\Gamma(a)=OD^{M_\Gamma}(\vec{A}, a)$. 
\end{proposition}
The proof of the observation above is straightforward and uses the fact that 
\begin{center}
$OD^{M_\Gamma}(\vec{A}, a)=\Sigma^2_1(\vec{A}, a)$.
\end{center}

Suppose now that $\Gamma$ is a vg-pointclass. We then say that $\vec{B}\subseteq \bR^{\omega}$ is a \textit{weakly} $\Gamma$-\textit{condensing sequence} if 
\begin{enumerate}
\item $B_0$ codes a sjs such that $M_\Gamma=L_\gg(B_0, \bR)$ and $\Gamma=(\Sigma^2_1(B_0))^{L_\gg(B_0, \bR)}$,
\item $B_1=\{ (x, y)\in\bR^2: y\in C_\Gamma(x)\}$, 
\item $B_2=B_1^c$,
\item $B_3$ is any $OD^{M_\Gamma}(B_0)$ set\footnote{We will need the freedom to include any $OD^{M_\Gamma}(B_0)$ set of reals into our condensing sequence.},
\item $(B_{2i+1}: i\in [2, \omega))\subseteq \Gamma$ is a scale on $B_1$,
 \item $(B_{2i}: i\in [2, \omega))\subseteq M_\Gamma$ is a scale on $B_2$,
 \item for every $i\in [2, \omega)$, $M_\Gamma\models ``B_{2i}$ is ordinal definable from $B_0"$.
\end{enumerate}

Suppose $\Gamma$ is a vg-pointclass and $M_\Gamma=L_\gg(\vec{A}, \bR)$ is its parent. Suppose $B\in M_\Gamma\cap \powerset(\bR)$ is $OD^{M_\Gamma}(\vec{A})$ and suppose $a\in HC$ is a transitive set. Consider the term relation $\tau^a_B$ consisting of pairs $(p, \sigma)$ such that 
\begin{enumerate}
\item $p\in Coll(\omega, a)$,
\item $\sigma\in C_\Gamma(a)$ is a standard $Coll(\omega, a)$-name for a real, and
\item for a co-meager many $g\subseteq Coll(\omega, a)$ (in the relevant topology) such that $p\in g$, $\sigma(g)\in B$\footnote{$\sigma(g)$ is the realization of $\sigma$.}.
\end{enumerate}
Then because $\tau^a_B$ is $OD^{M_\Gamma}(\vec{A}, a)$, $\tau^a_B\in C_\Gamma(C_\Gamma(a))$. Given $k\in \omega$, we let $\tau^a_{B, k}=\tau^{C^{k}_\Gamma}_{B, 0}$. Thus, for every $k\in \omega$, $\tau^a_{B, k}\in C^{k+2}_\Gamma(a)$.

We say $\vec{B}$ is a $\Gamma$-\textit{condensing sequence} if it is a weakly condensing sequence with the additional property that for any transitive sets $a, b, M\in HC$ such that
\begin{enumerate}
\item $a\in M$ and
\item there is an embedding $\pi: M\rightarrow_{\Sigma_1} C^\omega_\Gamma(b)$ such that $\pi(a)=b$ and for every $i, k\in \omega$, $\tau^b_{B_i, k}\in rng(\pi)$,
\end{enumerate}
$M=C^\omega(a)$ and for any $i, k\in \omega$, $\pi^{-1}(\tau^b_{B_i, k})=\tau^a_{B_i, k}$.

If $P$ is a $\Gamma$-Woodin and $B\in OD^{M_\Gamma}(\vec{A})$ then for $k\in \omega$ we let $\tau^P_{B, k}=\tau^{V_{\d^P}^P}_{B, k}$. 
\begin{definition} Suppose $\Gamma$ is a vg-pointclass and $M_\Gamma=L_\gg(\vec{A}, \bR)$ is its parent. Suppose $P$ is a $\Gamma$-Woodin and $\Sigma$ is an $\omega_1$-iteration strategy for $P$. Suppose $B\in M_\Gamma\cap \powerset(\bR)$ is $OD^{M_\Gamma}(\vec{A})$. 
\begin{enumerate}
\item We say $\Sigma$ is a $\Gamma$-fullness preserving strategy for $P$ if whenever $i:P\rightarrow Q$ is an iteration of $P$ via $\Sigma$, $Q$ is a $\Gamma$-Woodin. 
\item Given that $\Sigma$ is $\Gamma$-fullness preserving, we say $\Sigma$ respects $B$ if whenever $i:P\rightarrow Q$ is an iteration of $P$ via $\Sigma$, for every $k$, $i(\tau^P_{B, k})=\tau^Q_{B, k}$.  
\end{enumerate}
\end{definition} 

The following theorem, which probably is originally due to Woodin, is unfortunately unpublished. The discussion that follows the theorem might be illuminating. 

\begin{theorem}\label{coarse capturing} Assume $AD^+$ and suppose $\Gamma$ is a vg-pointclass. Let $M_\Gamma=L_\gg(\vec{A}, \bR)$\footnote{Here $\vec{A}$ is any sjs for which $M_\Gamma=L_\gg(\vec{A}, \bR)$.} be its parent and $A\in OD(\vec{A})^{M_\Gamma}$.  There is then a pair $(P, \Sigma)$ and a $\Gamma$-condensing sequence $\vec{B}$ such that
\begin{enumerate}
\item $P$ is a $\Gamma$-Woodin,
\item $\Sigma$ is a $\Gamma$-fullness preserving $\omega_1$-iteration strategy for $P$,
\item for each $i$, $\Sigma$ respects $B_i$,
\item for every $\Sigma$-iterate $Q$ of $P$, for every $i\in \omega$ and for every $Q$-generic $g\subseteq Coll(\omega, \d^Q)$, $\tau^Q_i(g)=Q[g]\cap B_i$,
\item for any tree $\T\in dom(\Sigma)$, $\Sigma(\T)=b$ if and only if either
\begin{enumerate}
\item $C_\Gamma(\M(\T))\models ``\d(\T)$ is not a Woodin cardinal" and $b$ is the unique well-founded cofinal branch $c$ of $\T$ such that $C_\Gamma(\M(\T))\in \M^\T_c$, or
\item $C_\Gamma(\M(\T))\models ``\d(\T)$ is a Woodin cardinal" and $b$ is the unique well-founded cofinal branch $c$ of $\T$ such that letting $Q=C^\omega_\Gamma(\M(\T))$, $\M^\T_c=Q$ and for every $i\in \omega$, $\pi^\T_c(\tau^P_{B_i})=\tau^Q_{B_i}$.
\end{enumerate}
\item $\Sigma$ respects $A$,
\end{enumerate}
Moreover, for any set $a\in HC$ there is $(P, \Sigma)$ as above such that $a\in P$. 
\end{theorem}

\begin{definition}\label{excellent pair} Suppose $\Gamma$ is a vg-pointclass. We say that $(P, \Sigma)$ is a $\Gamma$-excellent pair if for some $\Gamma$-condensing sequence $\vec{B}$, $(P, \Sigma)$ has properties 1-5 described in \rthm{coarse capturing} as witnessed by $\vec{B}$.
\end{definition}

The author believes that there is enough facts that have appeared in the literature that a motivated reader can put together a proof of \rthm{coarse capturing}. We now give a general outline of how the proof of \rthm{coarse capturing} might proceed.

First and foremost,  at the very heart of \rthm{coarse capturing} is Woodin's theorem that for a set of ordinals $T$ and for $T$-cone of $d$, $\omega_2^{L[T, d]}$ is a Woodin cardinal in $\H^{L[T, d]}(T)$ (see \cite[Theorem 5.4]{KoelWoodin}). The use of this theorem is as follows. Let $\Gamma$ be a vg-pointclass. Find a bigger vg-pointclass $\Gamma^*$ such that $\Gamma\subseteq \Delta_{\Gamma^*}$. Let $T$ be the tree of a $\Gamma^*$-scale on a $\Gamma^*$-universal set. Theorem 5.4 of \cite{KoelWoodin} shows that for any $x\in \bR$, we can find transitive $Q$ such that $x\in Q$ and letting $\nu=Ord\cap Q$, $L[T, Q]\models ``\nu$ is a Woodin cardinal". $\nu$ in question is $\omega_2^{L[T, z]}$ for some $z$ and $Q$ is $V_\nu^{\H^{L[T,z]}(T)}$. If we choose $z$ big enough then $C_\Gamma(Q)\in L[T, Q]$. 

The $\Gamma$-condensing sequence $\vec{B}$ can be produced via analyzing scales in $M_\Gamma$. It is a theorem of Martin that any set in $M_\Gamma$ has a semi-scale consisting of sets in $\Gamma$ (see \cite[Chapter 3]{Jackson}). We can then let $\vec{B}$ be any sequence such that for each $i\in \omega$, some subsequence $(B_{n_j}: j\in \omega)$ is a semi-scale on $B_i$. The proof of  \cite[Lemma 3.7]{MouseSetWoodin} shows that $\vec{B}$ is $\Gamma$-condensing. 

At this point, we can aim at obtaining a strategy for a $\Gamma$-Woodin. Let $\d$ be least such that $C_{\Gamma}(V_\d^Q)\models ``\d$ is a Woodin cardinal". Our choice of $z$ and $\Gamma^*$ guarantee that  $\d<\eta$. To see this, notice that $\{\tau^Q_{B_i, k}: i, k\in \omega\}\in L[T, Q]$ (consequence of choosing $\Gamma^*$ bigger than $\Gamma$). We now take a transitive below $\eta$ Skolem hull of $L[T, Q]$ such that the set $\{\tau^Q_{B_i, k}: i, k\in \omega\}$ is in it. Considering the collapse of the hull, we have $\pi: N\rightarrow L_\xi[T, Q]$. We now need to use the condensation properties of  $\{\tau^Q_{B_i, k}: i, k\in \omega\}$. These properties imply that if $\d=\cp(\pi)$ then $C_\Gamma(V_\d^N)\models ``\d$ is a Woodin cardinal" as $C_\Gamma(V_\d^N)\in N$. To gain more intuition on how exactly such condensation properties might work consider \cite[Lemma 3.7]{MouseSetWoodin}. There are such condensation lemmas at various parts of \cite{CMI} and \cite{ATHM}. This finishes our general outline of \rthm{coarse capturing}.

\rthm{coarse capturing} has an interesting consequence.

\begin{theorem}[Steel-Woodin]\label{measurable} Assume $AD^+$. Then all uncountable regular cardinals $<\Theta$ are measurable.
\end{theorem}

\rthm{measurable} can be proved via exactly the same method that is used to prove \cite[Theorem 8.27]{OIMT}. However, both theorems need a comparison theory which does not exist for $\Gamma$-Woodins (this is because $\Gamma$-Woodins are coarse models of a fragment of set theory not fine structural models). To actually implement the above mentioned proofs one needs to use fine structural models. \rthm{m_n exists} demonstrates how to obtain such models.

\section{$\Theta$-extensions of iteration strategies}

The following is the main theorem of this section. Assume $AD_{\mathbb{R}}$ and suppose $\Sigma$ is an $\omega_1$-iteration strategy for some $P$. Given a set $A$ such that 
\begin{enumerate}
\item $A\in L[A]$,
\item $L[A]\models ZFC$,
\item $A\in L_\b[A]$ for some $\b<\Theta$,
\end{enumerate}
we let $\b_A$ be the least $L[A]$-cardinal $\gg$ such that $A\in L_\gg[A]$ and $\mu_A$ be the $\omega_1$-supercompactness measure on $\powerset_{\omega_1}(L_{\b_A}[A])$\footnote{Recall that $AD_{\mathbb{R}}$ implies that $\omega_1$ is supercompact.}.  

We say $\Sigma^+$ is a $\Theta$-\textit{extension} of $\Sigma$ if for any $\T$ that is according to $\Sigma^+$ and is of length $<\Theta$, for a $\mu_{\T}$ measure one set of countable $X\prec L_\b[\T]$, letting $\T_X$ be the transitive collapse of $X$, $\T_X$ is according to $\Sigma$. 

Recall the definition of hull condensation from \cite[Definition 1.30]{ATHM}. Essentially hull condensation for a strategy $\Sigma$ is the statement that if $\T$ is according to $\Sigma$ and $\U$ is a Skolem hull of $\T$ then $\U$ is according to $\Sigma$. 

\begin{theorem}\label{main theorem 1} Assume $AD_{\mathbb{R}}$ and suppose $\Sigma$ is an $\omega_1$-iteration strategy for some $P$ with hull condensation. Then $\Sigma$ can be uniquely extended to a $\Theta$-strategy. Moreover, if $\Sigma$ is an $(\omega_1, \omega_1)$-iteration strategy with hull condensation then $\Sigma$ can be uniquely extended to a $(\Theta, \Theta)$-strategy. 
\end{theorem}
\begin{proof} We only prove the first part of the claim. 
Given $\T$ on $\M$ such that $lh(\T)<\Theta$, set\\\\
$Correct(\T)$: for a $\mu_{\T}$ measure one set of $X\prec L_\b[\T]$, letting $\T_X$ be the transitive collapse of $\T$, $\T_X$ is according to $\Sigma$\\\\

Suppose $Correct(\T)$ holds and $\T$ is of limit length.

\begin{lemma} There is a unique well-founded branch $b$ of $\T$ such that $Correct(\T^\frown\{b\})$ holds.
\end{lemma}
\begin{proof} Suppose $\cf(\T)>\omega$. It follows from \rthm{measurable} that $\cf(\T)$ is a measurable cardinal, and hence $\T$ has a unique branch $b$. 
We claim that $Correct(\T^\frown\{b\})$ holds. Let $\b=\b_{\T^\frown \{b\}}$ and fix a countable $X\prec  L_\b[\T^\frown\{ b\}]$. Let $\gg=\sup(X\cap lh(\T))$ and let $Y\prec  L_\b[\T^\frown\{ b\}]$ be countable and such that $X\cup\{\gg\}\subseteq Y$. Let $c_X: X\rightarrow N_X$ and $c_Y:Y\rightarrow N_Y$ be the transitive collapses of $X$ and $Y$. Let $c_{X, Y}: N_X\rightarrow N_Y$ be given by $c_{X, Y}(c_X(u))=c_Y(u)$ for $u\in X$. Because $Correct(\T)$ holds, we can find $Y$ as above such that $c_Y(\T)$ is according to $\Sigma$. Notice now that $c_{X, Y}\rest (c_X(\T^\frown \{b\}))$ witnesses that $c_X(\T^\frown \{b\})$ is a hull of $c_Y(\T)\rest c_Y(\gg+1)$. Therefore, $\Sigma(c_X(\T))=c_X(b)$. 

Assume now that $\cf(\T)=\omega$. Set $\b=\b_{\T}$ and $\mu=\mu_\T$. Given a countable $X\prec L_\b[\T]$, we let $c_X: X\rightarrow N_X$ be the transitive collapse of $X$ and $b_X=\Lambda(c_X(\T))$. Define $b\subseteq lh(\T)$ by setting: $\a\in b$ if and only if for a $\mu$-measure one many $X$,
\begin{enumerate}
\item $X\prec L_\b[\T]$,
\item $\a\in X$, and
\item $c_X(\a)\in b_X$.
\end{enumerate}

\begin{sublemma}\label{b is wf} $b$ is a cofinal well-founded branch of $\T$.
\end{sublemma}
\begin{proof} Assume $b$ is not cofinal. Let $\gg<lh(\T)$ be such that $b\subseteq \gg$. For each $X\prec L_\b[\T]$ such that $\gg\in X$, let $\gg_X=\min(b_X-c_X(\gg))$. Let $\a_X=c_X^{-1}(\gg_X)$. Then $\a_X\in X$. Because $\mu$ is normal, we have that for a some $\a<lh(\T)$, for $\mu$-measure one set of $X$, $\a_X=\a$. But then $\a\in b$. 

Towards showing that $b$ is a well-founded branch, notice that we have that for a $\mu$-measure one many $X$, $c_X[b\cap X]$ is cofinal in $b_X$. Indeed, let $(\a_i: i<\omega)$ be an increasing cofinal sequence in $b\cap lh(\T)$. For $i<\omega$, let $B_i\in \mu$ be such that for $X\in B_i$, 
\begin{enumerate}
\item $X\prec L_\b[\T]$, 
\item $(\a_i: i<\omega)\subseteq X$ and
\item  $c_X(\a_i)\in b_X$. 
\end{enumerate}
Let $B=\cap B_i$. Fix now $X\in B$. Then $c_X(\a_i)\in b_X$ for every $i<\omega$. Moreover, $(c_X(\a_i): i<\omega)$ is cofinal in $c_X(lh(\T))$.   

Fix now $B\in \mu$ such that for every $X\in B$, $c_X[b\cap X]$ is cofinal in $b_X$. Let $X\prec L_\b[\T^\frown \{b\}]$ be countable and such that $X\cap L_\b[\T]\in B$. Set $Z=X\cap L_\b[\T]$. Let $c:X\rightarrow N$ be the transitive collapse of $X$. It is enough to see that $c(b)$ is a well-founded branch of $c(\T)$. Notice that $c(\T)=c_Z(\T)$ and $c(b)=c_Z[b\cap Z]$. It then follows that $c(b)=b_Z$. Since $b_Z=\Lambda(c_Z(\T))$, we are done. 
\end{proof}

The proof of \rslem{b is wf} also verifies that $Correct(\T^\frown\{b\})$ holds. Indeed, the set $B$ witnesses it as shown by the last computation, namely that $c(b)=b_Z$.
\end{proof}

We now define $\Sigma^+$ by setting 
\begin{enumerate}
\item $dom(\Sigma^+)=\{\T:\T$ has a limit length $<\Theta$ and $Correct(\T)$ holds$\}$ and
\item $\Sigma^+(\T)=b$ if and only if $b$ is the unique branch of $\T$ such that $Correct(\T^\frown \{b\})$ holds.
\end{enumerate}
\end{proof}

\section{$\M_n^{\#, \Sigma}$}

Suppose $\Gamma$ is a vg-pointclass and $(P, \Sigma)$ is a $\Gamma$-excellent pair. For $n\in \omega$ and $x\in HC$ such that $L[x]\models ``x$ is well-ordered" (such $x$ are called swo\footnote{self-well-ordered}), we say that $\M_n^{\#, \Sigma}(x)$ exists if there is a sound active $\Sigma$-premouse $\M$ that projects to $\omega$, has $n$ Woodin cardinals  and is $\omega_1$-iterable (as a $\Sigma$-premouse\footnote{i.e., the iterates of $\M$ are also $\Sigma$-premice}). 

\begin{theorem}\label{m_n exists} Assume $AD^+$. Suppose $\Gamma$ is a vg-pointclass and $(P, \Sigma)$ is a $\Gamma$-excellent  pair. Then for every $n\in \omega$, $\M_n^{\#, \Sigma}$ exists. 
\end{theorem}
\begin{proof}  The proof is via induction on $n$. We first prove the inductive step.

\begin{sublemma} Assume for every $x\in HC$, $\M_n^{\#, \Sigma}(x)$ exists. Then for every swo $x$, $\M_{n+1}^{\#, \Sigma}(x)$ exists. 
\end{sublemma}
\begin{proof}
Towards a contradiction assume that $\M_{n+1}^{\#, \Sigma}(x)$ doesn't exist for some swo $x\in HC$.
Notice that the set $A_n=\{ (z, y)\in \bR^2:  y\ \text{codes}\ \M_n^{\#, \Sigma}(z)\}$ is projective in $\Sigma$. Let $\Pi$ be a vg-pointclass such that $\Gamma\subseteq \utilde{\Delta}_\Pi$.  The key fact that allows us to prove our claim is that all sets of reals that are projective in $Code(\Sigma)$ are in $\utilde{\Delta}_\Pi$. This immediately follows from the definition of vg-pointclass. Let $M_\Pi=L_{\nu}(D, \bR)$ be the parent of $\Pi$. We then have that $Code(\Sigma)\in (\utilde{\Delta}^2_1(D))^{M_\Pi}$. Because $(\utilde{\Delta}^2_1(D))^{M_\Pi}$ is closed under the real quantification, we have that all sets projective in $Code(\Sigma)$ are in $(\utilde{\Delta}^2_1(D))^{M_\Pi}$ and moreover, there is a set $U\in (\utilde{\Delta}^2_1(D))^{M_\Pi}$ that codes a sequence $(U_n: n<\omega)$ such that for each $n\in \omega$, $U_n$ is a universal $\Sigma^1_n(Code(\Sigma))$ set of reals. 

Let $(Q, \Lambda)$ be a $\Pi$-excellent pair such that $x\in Q$ and if $\vec{B}$ is $\Pi$-qsjs then $B_0$ codes the pair $(Code(\Sigma), A_n)$. Let $\d=\d^Q$. Notice that we have that $P\in HC^Q$. Let $\M$ be the fully backgrounded $L[\vec{E}][x]$-constructions relative to $\Sigma$ done inside $V_{\d^Q}^Q$. $\M$ inherits an $\omega_1$-strategy $\Psi$ from $\Lambda$ (see \cite[Chapter 12]{FSIT}). 

 Because $B_0$ codes $A_n$, we have that $\M_n^{\#, \Sigma}(V_\d^Q)\in Q$.Therefore, $\M_n^{\#, \Sigma}(V_\d^Q)\models ``\d$ is a Woodin cardinal". It follows that $\M_n^{\#, \Sigma}(\M)\models ``\d$ is a Woodin cardinal".\footnote{To see this, first because $\d$ is inaccessible, condensation implies that $\rho(\M_n^{\#, \Sigma}(\M))=\d$ as otherwise we could take a transitive below $\d$ hull $\N$ of $\M_n^{\#, \Sigma}(\M)$, and by condensation $\N\insegeq \M$. Next, fix $f:\d\rightarrow \d$ such that $f\in \M_n^{\#, \Sigma}(\M)$. Then following the proof of \cite[Theorem on page 115]{FSIT} find an extender $E$ such that $\pi_E^Q(f)(\cp(E))<strength(E)$ and $E\cap \M \in \M$. This $E$ witnesses Woodinness of $\d$ with respect to $f$.}. 
 
 However, because $\M$ is $\omega_1$-iterable as witnessed by $\Psi$, there is no active initial segement of $\M$ that reaches $n+1$ Woodin cardinals. Because $\M_{n}^{\#, \Sigma}(\M)$ has $n+1$ many Woodins and is active, we have that $\rho(\M_{n}^{\#, \Sigma}(\M))=\omega$. However, as we argued in the footnote below, $\rho(\M_n^{\#, \Sigma}(\M))=\d$. 
\end{proof}

To finish the proof of the lemma, we need to show that $\M_0^{\#, \Sigma}(x)$ exists for every $x$. The proof of the sublemma above shows how to do this. All we need to do is get any $(Q, \Lambda)$ as in the proof of sublemma and consider $\M$ defined the same way as in the proof of the sublemma. Because there are measurable cardinals in $Q$, we must have that $\M$ has measurable cardinals. It follows that $\M_0^{\#, \Sigma}(x)$ exists.
\end{proof}

\section{The internal theory of $\M_1^{\#, \Sigma}$} 

Fix a vg-pointclass $\Gamma$ and a $\Gamma$-excellent pair $(P, \Sigma)$. Fix $\vec{B}$ witnessing the $\Gamma$-excellence of $(\P, \Sigma)$ and let $u$ be a real coding the sequence $(\tau^P_{B_i}: i<\omega)$. Let $\M=\M_1^{\#, \Sigma}(u)$ and let $\Lambda$ be the $\omega_1$-iteration strategy of $\M$. The following is the main theorem of this section.

\begin{theorem}\label{internal uB representation} Let $\d$ be the Woodin cardinal of $\M$. There are trees $(T, S)\in \M$ on $\omega\times (\d^+)^\M$ such that $\M\models ``(T, S)$ are $\d$-complementing" and whenever $i: \M\rightarrow \N$ is an iteration according to $\Lambda$ and $g\subseteq Coll(\omega, i(\d))$ is $\N$-generic, $\N[g]\cap p[i(T)]=Code(\Sigma)\cap \N[g]$. 
\end{theorem}
\begin{proof}
The proof uses the idea of \textit{generic genericity} iterations. Fixing a cardinal $\eta$ of $\M$, there is a tree $\T_\eta\in \M$ on $P$ that is according to $\Sigma$ with last model $P_\eta$ such that for any real $x$ that is generic over $\M$ for a poset in $\M|\eta$, $x$ is generic over $\mathbb{W}_{\d^{P_\eta}}^{P_\eta}$\footnote{$\mathbb{W}_\eta^S$ is the extender algebra of $S$ associated to $\eta$, a Woodin cardinal of $S$.}. To obtain $\T_\eta$ first iterate the least measure of $P$, $\eta+1$ times and then follow the usual construction of genericity iterations with a slight modification; namely that at a successor step of the genericity iteration we pick the least extender $E$ such that there is a condition $p\in Coll(\omega, \eta)$ that forces that some real $x$ violates some axiom generated by $E$.  

Let $\T=\T_{\d^+}$ and $Q$ be the last model of $\T$. Working in $\M$, we want to find a formula $\phi$, a formula $\psi$ and a parameter $r\in \M|(\d^{+2})^\M$ such that\\\\
$\dagger(\phi, \psi, r):$ for a club of countable $X\prec \M|(\d^{+2})^\M$ letting $\pi_X: \N_X\rightarrow \M|(\d^+)^\M$ be the transitive collapse of $X$ and setting $\d_X=\pi_X^{-1}(\d)$ and $r_X=\pi^{-1}_X(r)$, for any $\mathbb{P}\in \N_X|(\d_X^+)^{\N_X}$ whenever $g\subseteq \mathbb{P}$ is $\N_X$-generic the following conditions hold:
\begin{enumerate}
\item For every $y\in \mathbb{R}\cap \N_X[g]$, $\N_X[g]\models \phi[y, r_X]$ if and only if $y$ codes a tree $\U\in dom(\Sigma)$.
\item For every $(y, z)\in \mathbb{R}^2\cap \N_X[g]$, $\N_X[g]\models \psi[y, z, r_X]$ if and only if $y$ codes a tree $\U\in dom(\Sigma)$ and $z$ codes $\Sigma(\U)$.
\item $\N_X[g]\models \forall y\in \bR(\phi(y, r_X)\rightarrow \exists z\in \bR \psi(y, z, r_X))$.
\end{enumerate} 
We will take $r=\T$. We describe $\phi$ and $\psi$ by first describing another formula $\sigma$. Notice that the language of $\M$ already has names  for $\vec{E}^\M$, $P$, $u$ and $\Sigma$. Below we describe $\sigma$, $\phi$ and $\psi$.\\\\
\textbf{The description of $\sigma$:} We let $\sigma(y, r)$ be the formula expressing the following statement: $y\in \bR$ codes a tree $\U$ on $P$ such that for every limit $\a<lh(\U)$ letting $b=[0, \a)_\U$, $b$ is the unique branch $c$ of $\U\rest \a$ such that either
\begin{enumerate}
\item $C_\Gamma(\M(\U\rest \a))\models ``\d(\U\rest \a)$ is not a Woodin cardinal", $\M^{\U\rest \a}_c$ is well-founded and $C_\Gamma(\M(\U\rest \a))\in \M^{\U\rest \a}_c$ or
\item $C_\Gamma(\M(\U\rest \a))\models ``\d(\U\rest \a)$ is a Woodin cardinal", $\M^{\U\rest \a}_c$ is well-founded, and letting $R=C_\Gamma(C_\Gamma(\M(\U\rest \a))$, $R=\M^{\U\rest \a}_c$ and $\pi^{\U\rest \a}_c(\tau^P_{B_i})=\tau^R_{B_i}$ for every $i$.\\
\end{enumerate}
\textbf{The description of $\phi$:} We let $\phi(y, r)$ be the formula expressing the following statement: $\sigma(y, r)$ and if $\U$ is the tree coded by $y$ then $\U$ has a limit length. \\\\
\textbf{The description of $\psi$:} We let $\psi(y, z, r)$ be the formula expressing the following statement: $\phi(y, r)$ and if $\U$ is the tree coded by $y$ then $z\in \bR$ codes a cofinal well-founded branch $b$ of $\U$ such that letting $w$ be any real coding $\U^\frown\{b\}$, $\sigma(w, r)$ holds. \\\\
Clearly the definition of $\sigma$ is vague as it refers to $C_\Gamma$ operator and the sequence of terms $(\tau^R_{B_i}: i\in \omega)$. We now describe  formulas $\sigma_0$ and $\sigma_1$ that define the functions $x\rightarrow C_\Gamma(x)$ and $R\rightarrow (\tau^R_{B_i}: i\in \omega)$. The appropriate structures for evaluating $\sigma_0$ and $\sigma_1$ are again small generic extensions of Skolem hulls of $(\M|(\d^{+2})^\M, \T, \in)$.\\\\
\textbf{The description of $\sigma_0$:} We let $\sigma_0(y, z, \T)$ be the formula expressing the following statement: $(y, z)\in \bR^2$ and $z$ codes the set of all reals $w$ such that 
\begin{center}
$Q[y, w]\models \forces_{Coll(\omega, \d^Q)} (\check{y}, \check{w}) \in \pi^\T(\tau_{B_0}^P)$.\\ 
\end{center}
\textbf{The description of $\sigma_1$:} We let $\sigma_1(y, i, z, w, \T)$ be the formula expressing the following statement: $i\in \omega$, $(y, z, w)\in \bR^3$, $\sigma_1(y, z, \T)$ holds, and letting $a$ be the set of reals coded by $z$, $w$ codes the set of those reals $t\in a$ such that  
\begin{center}
$Q[t]\models \forces_{Coll(\omega, \d^Q)} t \in \pi^\T(\tau_{B_i}^P)$. 
\end{center}
 We then have the following two lemmas.
 \begin{lemma}\label{cgamma ub} Suppose $X\prec (\M|(\d^{+2})^\M, \T, \in)$ is countable. Let $\pi_X: \N_X\rightarrow \M|(\d^+)^\M$ be the transitive collapse of $X$ and set $\d_X=\pi_X^{-1}(\d)$, $\T_X=\pi^{-1}_X(\T)$ and $Q_X=\pi^{-1}_X(Q)$. Let $\mathbb{P}\in \N_X|(\d_X^+)^{\N_X}$ be a poset and  $g\subseteq \mathbb{P}$ be $\N_X$-generic. Suppose further that $a\in HC^{Q_X[g]}$ is an swo. 
 \begin{enumerate}
 \item Assume $a\subseteq \omega$. Working in $\N_X[g]$, let $b$ the set of of all those reals that are coded by any $z$ such that $\N_X[g]\models \sigma_0[a, z, \T_X]$. Then $b=C_\Gamma(a)$. 
 \item Assume $a$ is a transitive model of some fragment of $ZFC$. Working in $\N_X[g]$, let $b$ the set of of all those sets $c^*$ that are coded by some  $c\subseteq a$ with the property that for a comeager many $h\subseteq Coll(\omega, a)$, for any $z$ such that $\N_X[g]\models \sigma_0[a_h, z, \T_X]$, $ c_h$ is one of the reals coded by $z$. Then $b=C_\Gamma(a)$. 
 \end{enumerate}
 \end{lemma}
 \begin{proof} Clause 1 is easy and follows from the fact that $\vec{B}$ is a $\Gamma$-condensing sequence. Decoding $\sigma_0$ shows that what we have is that $b$ is the set of all those sets $c^*$ that are coded by some $c\subseteq a$ with the property that $\N_X[g]\models ``$for a comeager many $h\subseteq Coll(\omega, a)$, $c_h\in C_\Gamma(a_h)$". This later fact is absolute between $\N_x[g]$ and $V$. Hence, $b=C_\Gamma(a_h)$. 
 \end{proof}
The proof of the next lemma is similar and we leave it to the reader.
 \begin{lemma}\label{internal capturing of b} Suppose $X\prec (\M|(\d^{+2})^\M, \T, \in)$ is countable. Let $\pi_X: \N_X\rightarrow \M|(\d^+)^\M$ be the transitive collapse of $X$ and set $\d_X=\pi_X^{-1}(\d)$, $\T_X=\pi^{-1}_X(\T)$ and $Q_X=\pi^{-1}_X(Q)$. Let $\mathbb{P}\in \N_X|(\d_X^+)^{\N_X}$ be a poset and  $g\subseteq \mathbb{P}$ be $\N_X$-generic. Suppose further that $a\in HC^{Q_X[g]}$ is an swo and $b=C_\Gamma(a)$. Fix $i\in \omega$ and let $\tau\in b^{Coll(\omega, a)}$ be a name consisting of pairs $(\theta, p)$ such that in $\N_X[g]$, 
 \begin{enumerate}
 \item $\theta$ is a standard name for a real and 
 \item for comeager many $h\subseteq Coll(\omega, a)$, for any $z$ coding $b$, $\N_X[g]\models \sigma_1(i, a_h, z,\theta(h))$ holds. 
 \end{enumerate}
  Then $\tau=\tau^b_{B_i}$ and $\tau\in C_\Gamma(b)$. \\
 \end{lemma}
 
 We now resume the proof of \rthm{internal uB representation}. First we verify that $\dagger(\phi, \psi, \T)$ holds. To see this, fix $X\prec (\M|(\d^{+2})^\M, \T, \in)$ and let $\pi_X:\N_X\rightarrow (\M|(\d^{+2})^\M, \T, \in)$ be the transitive collapse of $X$. We will use subscript $X$ to denote $\pi_X$-preimages of objects in $X$. Let $\mathbb{P}\in \N_X|(\d_X^+)^{\N_X}$ be a poset and $g\subseteq \mathbb{P}$ be $\N_X$-generic in $\M$.The following is the main claim towards the proof of \rlem{internal uB representation}.\\
 
 \textit{Claim 1.} Suppose $\U\in HC^{\N_X[g]}\cap dom(\Sigma)$. Then $\Sigma(\U)\in \N_X[g]$.\\\\
 \begin{proof} Let $\Sigma(\U)=b$. We then know that $b$ is the unique well-founded branch $c$ of $\U$ such that for every $i<\omega$, $\pi^\U_c(\tau^P_{B_i})=\tau^Q_{B_i}$ where $Q=C_\Gamma(C_\Gamma(\M(\U))$. \rlem{cgamma ub} and \rlem{internal capturing of b} imply that $Q, b\in \N_X[g]$.
 \end{proof}
Fix now $y\in \bR^{\N_X[g]}$. \\

\textit{Claim 2.} Suppose $y$ codes a tree $\U\in dom(\Sigma)$ and $\N_X[g]\models \phi[y, \T]$. Then letting $b=\Sigma(\U)$, $b\in \N_X[g]$ is the unique well-founded branch $c$ of $\U$ such that whenever $z\in \N_X[g]\cap \bR$ codes $\U^\frown\{ c\}$, $\N_X\models \phi[z, \T]$. \\\\
\begin{proof} We only need to verify that $b$ satisfies clauses 1 and 2 of $\sigma[z, \T]$. But this immediately follows from the fact that $\Sigma(\U)=b$,  \rlem{cgamma ub} and \rlem{internal capturing of b}. 
\end{proof}

\textit{Claim 3.} $\N_X[g]\models \phi[y, \T]$ if and only if $y$ codes a tree $\U\in dom(\Sigma)$.\\\\
\begin{proof} That $y$ codes a tree on $P$ is absolute between $\N_X[g]$ and $\M$. We then assume that $y$ indeed codes a tree $\U$ of limit length. 

Assume first that $\N_X[g]\models \phi[y, \T]$. Fix a limit $\a<lh(\U)$ and assume that $\U\rest \a\in dom(\Sigma)$. Let $b=[0, \a)_\U$. It follows from \rlem{cgamma ub} and \rlem{internal capturing of b} that $\Sigma(\U\rest \a)=b$.

Suppose next that $\U\in dom(\Sigma)$. Let $\a<lh(\U)$ be a limit ordinal and let $b=[0, \a)_\U$. Because $b=\Sigma(\U\rest \a)$, we have that the two clauses of formula $\sigma$ are satisfied in $\M$. \rlem{cgamma ub} and \rlem{internal capturing of b} imply that the two clauses of $\sigma$ are also satisfied in $\N_X[g]$. 
\end{proof}

The next claim is an easy corollary of Claim 1-3.\\

\textit{Claim 4.} Suppose $\N_X[g]\models \phi[y, \T]$. Then there is a real $z\in \N_X[g]$ such that $\N_X\models \psi[y, z, \T]$.\\\\
We now have that Claim 1-4 imply that $\dagger(\phi, \psi, \T)$ holds. It follows that we can find $(T, S)\in \M|(\d^{+2})^M$ such that 
\begin{enumerate}
\item $\M\models ``(T, S)$ are $(\d^+)^\M$-complementing", and
\item for every $<(\d^+)^\M$-generic $g$, in $\M[g]$, $(y, z)\in p[T]$ if and only if there is a countable $X\prec (\M[g], \T, \in)$ such that $(y, z)\in X$ and letting \begin{center}
$\pi_X:\N_X\rightarrow (\M[g], \T, \in)$
\end{center}
 be the transitive collapse of $X$, $\N_X\models \psi[y, z, \pi_X^{-1}(\T)]$.
\end{enumerate}
Repeating the proofs of Claim 1-4 and using elementarily we see that whenever $i:\M\rightarrow \N$ is an iteration via $\Lambda$,\\\\
(*) for every $<i((\d^+)^\M)$-generic $g$, $p[i(T)]\cap \N[g]=Code(\Sigma)\cap \N[g]$.\\\\
Clearly (*) finishes the proof of \rlem{internal uB representation}.
\end{proof}

\section{All sets are $\Theta$-uB}

In this section, we prove \rthm{main theorem}. Assume $AD_{\mathbb{R}}$. We will describe two approaches and complete only one leaving the other to the reader. Both approaches are the same in the sense that they use the same ideas. However, approach 2 uses the theorems and lemmas that we have already proved. Nevertheless, Approach 1 is easier and doesn't use the material on the existence of iterable models with Woodin cardinals. \\

\textbf{Approach 1}\\

Let $\mu$ be the supercompactness measure on $\powerset_{\omega_1}(\bR)$. Fix a set of reals $A$ and find a vg-pointclass $\Gamma_0$ such that it has a parent $M_{\Gamma_0}=L_{\gg}(\vec{A}, \bR)$ with $A=A_0$. Let $\Gamma$ be a vg-poiintclass such that $\Gamma_0\subseteq \Delta_{\Gamma}$ and $\Gamma$ has a parent of the form $L_\xi(\vec{C}, \bR)$ such that $C_0$ codes the sequence $\vec{B}$.  Let $\vec{B}$ be a $\Gamma$-condensing sequence. Thus $A$ is coded into $\vec{B}$.
 
 Let $T$ be a tree of a $\Gamma$-scale on a $\Gamma$-universal set. We say $Q\in HC$ is $\Gamma$-\textit{full} if letting $o(Q)=Ord\cap Q$, $Q=V_{o(Q)}^{L[T, Q]}$. 

Suppose now that $Q$ is $\Gamma$-full and is a model of $ZF$. We then have that\\\\
(1) for every set $x\in Q$, $C^\omega_\Gamma(x)\in Q$,\\
(2) for every set $x\in Q$, $\{ \tau^x_{B_i, k}: i, k\in \omega\}\in Q$.\\\\
Give a transitive set $a\in HC$, let $\xi_a<\omega_1$ be least such that $L_{\xi_a}(T, a)\models ZF$ and set $Q_a=L_{\xi_a}(T, a)$. We have that $Q_a$ is $\Gamma$-full (see \rthm{capturing by trees}). Using $\mu$, we can extend the function $a\rightarrow Q_a$ to all transitive sets in $L_{\Theta}(\powerset(\bR))$ by repeating the proof of \rthm{main theorem 1}. By repeating the proof of \rthm{internal uB representation}, we can show that
\begin{lemma}\label{capturing in qa} for every transitive $a\in HC$, there are trees $(T_a, S_a)\in Q_a$ such that letting $\nu_a=\Theta^{Q_a}$ and $\l_a=(\nu_a^+)^{Q_a}$,
\begin{enumerate}
\item $Q_a\models ``(T_a, S_a)$ are $\nu_a$-complementing",
\item for any $\mathbb{P}\in Q_a$ such that $Q_a\models ``\mathbb{P}$ is a surjective image of $\bR^{Q_a}"$ and for any $Q_a$-generic $g\subseteq \mathbb{P}$,
\begin{center}
$Q_a[g]\cap p[T_a]=C_0\cap Q_a[g]$
\end{center} 
\end{enumerate}
\end{lemma}

Let now $Q^+: L_{\Theta}(\powerset(\bR))\rightarrow L_\Theta(\powerset(\bR))$ be the extension of $Q$ to all transitive sets in $L_{\Theta}(\powerset(\bR))$. The method extending functions defined on $HC$ to $L_{\Theta}(\powerset(\bR))$ given in the proof of \rlem{main theorem 1} guarantees that\\\\
(3) for $a$ a transitive set in $L_{\Theta}(\powerset(\bR))$, letting $\mu^*$ be the $\omega_1$-supercompactness measure on $\powerset_{\omega_1}(Q^+_a)$ obtained as a projection of $\mu$, for $\mu^*$-measure one many $X\prec (Q^+_a, \in)$ letting $\tau_X: a\rightarrow N_X$ be the transitive collapse of $X$, $N_X=Q_{\tau_X(a)}$.\\\\
It follows from \rlem{capturing in qa} that for every transitive $a\in L_{\Theta}(\powerset(\bR))$, there are trees $(T^+_a, S^+_a)\in Q^+_a$ such that \\\\
(4) letting $\mu^*$ be the $\omega_1$-supercompactness measure on $\powerset_{\omega_1}(Q^+_a)$ obtained as a projection of $\mu$, for $\mu^*$-measure one many $X\prec (Q^+_a, \in)$ letting $\tau_X: a\rightarrow N_X$ be the transitive collapse of $X$, $N_X=Q_{\tau_X(a)}$ and $\tau_X(T^+_a, S_a^+)=(T_a, S_a)$. \\\\
Combining (4) with \rlem{capturing in qa}, we get that\\\\
(5) for any transitive $a\in L_\Theta(\mathbb{R})$ model of $ZF$ that contains $\bR$, for any $\mathbb{P}\in Q^+_a$ such that $Q^+_a\models ``\mathbb{P}$ is a surjective image of $\bR$" and for any $Q^+_a$-generic $g\subseteq Coll(\omega, \bR)$, $(p[T])^{Q_a^+[g]}=(p[S]^{Q^+_a[g]})^c$ and $p[T]=C_0$.\\\\
Clearly (5) implies that $C_0$  is universally Baire. Hence, $A$ is also universally Baire.\\\\

\textbf{Approach 2}\\

As we said above, approach 2 is very similar. Let $\mu$ be the supercompactness measure on $\powerset_{\omega_1}(\bR)$. Fix a set of reals $A$ and find a vg-pointclass $\Gamma$ such that it has a parent $M_{\Gamma}=L_{\gg}(\vec{A}, \bR)$ with $A=A_0$. Let $\vec{B}$ be a $\Gamma$-condensing sequence such that $B_0$ codes $\vec{A}$ and let $(P, \Sigma)$ be a $\Gamma$-excellent pair as witnessed by $\vec{B}$.

Given a countable $\sigma\subseteq \bR$ and $C\subseteq \bR$, let $\N(C, \sigma)=\M_1^{\#, \Sigma}((C\cap \sigma, \sigma)^\#)$. Recall that under $AD_\bR$ every set of reals has a sharp\footnote{$C^\#$ can be thought as the minimal active mouse over $C, \bR$ that is $\omega_1$-iterable}. 

\begin{lemma}\label{extension} Let $F:\powerset(\bR)\times \powerset_{\omega_1}(\bR)\rightarrow HC$ be given by $F(C, \sigma)=\N(C, \sigma)$. Then there is a total $F^+: \powerset(\bR) \rightarrow L_{\Theta}(\powerset(\bR))$ such that for every $C\in dom(F^+)$, for a club of $X\prec (F^+(C), C, \in)$, letting $\pi_X: (F^+(C), C)\rightarrow (N^+_C, N_C)$ be the transitive collapse of $X$, $N^+_C=\N(C, \pi_X[\bR])$.
\end{lemma}
\begin{proof}
Given $C\subseteq \bR$, let $\mu_C$ be the supercompactness measure on $\powerset_{\omega_1}(C^{\#})$. Set $F^+(C)=(\Pi \N(C, C\cap \pi_X[\bR]))/ \mu_C$ where the product ranges over countable $X\prec C^\#$ and $\pi_X: C^\#\rightarrow N_X$ is the transitive collapse of $X$. 

We want to show that $F^+$ has the desired properties. Fix $C\subseteq \bR$. Suppose now that $Y\prec (F^+(C), C, \in)$ is countable. Let $(M^+_Y, M_Y, \in)$ be the transitive collapse of $Y$. We can then fix $(f_i: i<\omega)$ such that for every $i$, $f_i: \powerset_{\omega_1}(C^\#)\rightarrow HC$ is such that $\{[f_i]_{\mu}: i<\omega\}=Y$. There is then a $\mu_C$-measure one set $S_0$ of $X\in \powerset_{\omega_1}(C^\#)$ such that for each $i$, $f_i(X)\in \N(C, \pi_X[\bR])$ where $\pi_X: X\rightarrow N_X$ is the transitive collapse of $X$. It follows that for any such $X$, there is an embedding $\sigma_X: N^+_Y \rightarrow \N(C, \pi_X[\bR])$ given by $\sigma_X(\pi_Y(a))=f_i(X)$ where $a=[f_i]_{\mu_C}$. By shrinking $S_0$ further to some $S\in \mu_C$ we can make sure that $\sigma_X$ is $\Sigma_1$-elementary for every $X\in S$. It then follows that $M^+_Y=\N(C, \pi_Y[\bR])$
\end{proof}

For each $C\subseteq \bR$, let $(T_C, S_C)\in F^+(C)$ be the trees described in \rlem{internal uB representation}. For each $C$ and $X\prec  C^\#$ we let $f(X)=(T, S)$ where $(T, S)\in \N(C, \pi_X[\bR])$ are the trees described in the proof of \rthm{internal uB representation} and $\pi_X: X\rightarrow N_X$ is the transitive collapse of $X$. Let then $(T_C, S_C)=[f]_{\mu_C}$. It follows from the definition of $F^+(C)$ that\\\\
(1) For every poset $\mathbb{P}\in F^+(C)$ such that $F^+(C)\models ``\mathbb{P}$ is a surjective image of $\bR$", whenever $g\subseteq \bR$ is $F^+(C)$-generic, $F^+(C)[g]\models ``p[T_C]=(p[S_C])^c$.  \\\\
Notice that if $C\subseteq \bR$ is such that $C_0\in C^\#$\footnote{This could happen if for instance, $C$ is Wadge reducible to $C$.} then for a $\mu_C$-measure one many $X\prec C^\#$, $Code(\Sigma)\in X$ and letting $\pi_X: X\rightarrow N_X$ be the transitive collapse of $X$, 
\begin{center}
$(p[T])^{\N(C, \pi_X[\bR])}=Code(\Sigma)\cap \N(C, \pi_X[\bR])$.
\end{center}
It follows that \\\\
(2) for every $C\subseteq \bR$, $p[T_C]=Code(\Sigma)$.\\\\
Thus,\\\\
(3) for every $C$, $(T_C, S_C)$ is a universally Baire representation of $Code(\Sigma)$ for posets $\mathbb{P}\in C^\#$ with the property that $C^\#\models ``\mathbb{P}$ is a surjective image of $\bR"$.\\\\
As for every poset $\mathbb{P}$ that is a surjective image of a reals there is in some $C^\#$ such that $\mathbb{P}\in C^\#$ and $C^\#\models ``\mathbb{P}$ is a surjective image of $\bR"$, it follows from (3) that\\\\
(4) $Code(\Sigma)$ is $\Theta$-uB.\\\\
To conclude that $A$ is also $\Theta$-uB we use the fact that $\Sigma$-respects $A$. More precisely, the set $U$ consisting of $(x, y)\in \bR^2$ such that 
\begin{enumerate}
\item $x$ codes a tree $\T$ on $P$ with last model $Q$, 
\item $y$ is generic over the extender algebra of $Q$ at $\pi^\T(\d)$ and 
\item $Q[y]\models \forces_{Coll(\omega, i(\d^P))} \check{y}\in \pi^\T(\tau^P_{A, 0})$,
\end{enumerate}
is $\Pi^1_1(u)$ for some real $u$. It is the easy to extract a $\Theta$-uB representation for $A$ from one for $Code(\Sigma)$ and $U$. This finishes the proof of \rthm{main theorem}.

We finish this section with an application. The authors of \cite{AllSetsUB} produced models of determinacy that satisfy ``All sets of reals are universally Baire".  They worked under large cardinal assumptions and the proof used the proof of the derived model theorem. Here we use determinacy. First we record the following easy corollary to our construction of the universally Baire representation of iteration strategies.

\begin{corollary}\label{cor 1} Assume $AD_{\mathbb{R}}$ and suppose $\Gamma$ is a vg-pointclass. Let $(P, \Sigma)$ be a $\Gamma$-excellent pair. Then $Code(\Sigma)$ is $\Theta$-uB as witnessed by a pair of trees $(T, S)$ that is ordinal definable from $Code(\Sigma)$ and a real.
\end{corollary}

Here is our second corollary. 

\begin{corollary}\label{a model of all sets are ub} Assume $AD_{\mathbb{R}}$. Let $(\theta_\a: \a\leq \Omega)$ be the Solovay sequence, and fixing $\a<\Omega$, set $\Gamma=\{A: w(A)<\a\}$. Then for all $\b\in [\a+1, \Omega)$ such that $\H\models ``\theta_\b$ is regular", $V_{\theta_{\b}}^{\H(\Gamma)}\models ZF+AD^++$``All sets are universally Baire". 
\end{corollary}
\begin{proof}
It is known that if $\H\models ``\theta_\b$ is regular" then $\H(\Gamma_\b)\models ``\theta_\b$ is regular" where  $\Gamma_\b=\{A\subseteq\bR: w(A)<\theta_\b\}$ (see for instance \cite[Theorem 2.3]{SquareAim}). It follows that for each $\b$ as in the statement of the corollary, $V_{\theta_{\b}}^{\H(\Gamma)}\models ZF+AD^+$. That $V_{\theta_{\b}}^{\H(\Gamma)}\models$`` All sets are universally Baire" follows from \rcor{cor 1} and the short argument given after bullet point (4) in the proof of \rthm{main theorem}.
\end{proof}

\bibliographystyle{plain}
\bibliography{ADRUB.bib}
\end{document}